\UseRawInputEncoding

\documentstyle[12pt,twoside, epsf]{article}

\let \ttorg \tt \def \tt{\ttorg \obeyspaces}

\begin{document}

\date{}


\title{\bf Calculus, Gauge Theory and Noncommutative Worlds}

\author{Louis H. Kauffman \\
  Department of Mathematics, Statistics and Computer Science \\
  University of Illinois at Chicago \\
  851 South Morgan Street\\
  Chicago, IL, 60607-7045}

\maketitle
  
\thispagestyle{empty}

\noindent {\bf Abstract.}  This paper shows how gauge theoretic structures arise in a noncommutative calculus where the derivations are generated by commutators. These patterns include Hamilton's equations, the structure of the Levi-Civita connection, and generalizations of electromagnetism that are related to gauge theory and with the early work of Hermann Weyl. The territory here explored is self-contained mathematically. It is elementary, algebraic and subject
to possible generalizations that are discussed in the body of the paper.
\\

\noindent {\bf Keywords.} derivative, derivation, commutator, curvature, parallel translation, covariant derivative, connection, Levi-Civita connection, gauge connection, holonomy,
discrete derivative, quantum theory, quantum gravity\\

\noindent{\bf AMS Classification.} 53Z05.\\

\section{Introduction to noncommutative Worlds}
This paper is an exploration of calculus in a noncommutative framework and its relationships with physics. It is well known that quantum mechanics can be formulated in such a framework. Here we begin with classical physics and show that it is illuminated by thinking about a context of non-commutativity. It is of interest to ask the relationship of such work with the Noncommutative Geometry of  Alain Connes \cite{Connes}. We give a concise description of the Connes  approach in Section 2 of this paper. There is a close
conceptual relationship of this work with  Connes  since both approaches represent calculus algebraically. The reader familiar with the work of Connes may find Section 2 helpful in orienting to the present paper. The contents of this paper are self-contained and elementary. 
We intend to make contact with Connes' Geometry in subsequent work.\\

When we say that this paper is elementary we mean it. There is one key idea in the present paper. We are given an algebra ${\cal N}$ that is associative and noncommutative.
We note that in such an algebra, if $J$ is a chosen element of the algebra and we define $\nabla_{J}: {\cal N} \longrightarrow {\cal N}$ by the equation $\nabla_{J} F = FJ - JF = [F,J],$ then $\nabla_{J}$ is a {\it derivation} on {\cal N} in the sense that $\nabla_{J}(FG) = \nabla_{J}(F)(G) + F\nabla_{J}(G).$ Thus $\nabla_{J}$ satisfies the Leibniz rule for differentiation. This means that we can mimic differential calculus in the context of the {\it noncommutative world} of the algebra ${\cal N}.$ In this paper we will discuss this algebraic version of calculus in noncommutative worlds. We shall refer to the usual contexts of differential calculus with spaces, topologies, and limits as {\it standard worlds}.\\

The relationships between standard worlds and noncommutative worlds are only partially explored. Quantum mechanics begins such an exploration when it follows the Dirac dictum
and replaces the Poisson brackets of Hamiltonian mechanics by commutators so that the quantum evolution of a wavefunction $\psi$ is given by the equation 
\begin{equation}
i \hbar \dot{\psi} = [\psi, \hat{H}]
\end{equation}
where $\hat{H}$ is the Hamiltonian operator for the quantum system.
In this case the operator $\hat{H}$  has an Hermitian representation on an appropriate Hilbert space and the time evolution is unitary.
There is more given structure than in  our abstract noncommutative world, and one does not attempt to do all the calculus by using commutators.
Nevertheless, the fact that quantization is performed by replacing standard calculus (the Poisson bracket) by commutators is in back of our motivation to 
explore the noncommutative world. Another motivation is in the Feynman-Dyson derivation of electromagnetism from commutator equations \cite{Dyson}.
In that derivation, the authors essentially work in a noncommutative world. We will have more to say about the Feynman-Dyson derivation in Section 8.\\

There is no continuum differential calculus in this paper. There are no topological spaces. There are no bundles, no tangent spaces, no cotangent spaces,
no spaces at all. The contents of this paper are entirely algebraic. It is our intent to put flesh on these bones, but that work will be done elsewhere. Here we are exploring the consequences of 
derivations defined by commutators in an abstract algebra. What is remarkable is that many {\it patterns} of physics and gauge theory arise in this algebraic context. It is our intent to explore and
exhibit these patterns.\\

Calculus was originally formulated in a commutative framework by Newton, Leibniz  and their successors.
Quantum mechanics brought formulations of physical theory  linked with non-commutativity.
Heisenberg's quantum theory is based on quantities that do not commute with one another. These quantities obey specific identities such as the commutator equation for position $Q$ and momentum $P:$
\begin{equation}
QP - PQ = \hbar i.
\end{equation}
Schr\"{o}dinger formulated quantum mechanics via partial differential equations and showed that the operators
$Q=x$ and $P = -i \hbar \partial/\partial x$ obey the Heisenberg relations.
Dirac found a key to quantization via the replacement of the Poisson bracket (of Hamiltonian mechanics) with the commutators of quantum operators.
Curvatures in differential geometry and general relativity are seen, through work of Weyl and others to correspond to differences in parallel translation. This corresponds to the commutators of covariant derivatives. 
Gauge theory began, with Hermann Weyl \cite{Weyl} as a generalization of differential geometry. In Weyl's theory, lengths as well as angles are dependent upon the choice of paths. Weyl saw how to incorporate electromagnetism into general relativity using his generalizations of differential geometry.  Initial difficulties of interpretation arose in the context of Weyl's theory for general relativity, but his ideas  were adopted for quantum mechanics and became a basis for understanding nuclear forces.\\

We begin by formulating calculus in noncommutative domains. Discrete calculus is a motivation for these constructions. We show how to embed discrete calculus in a noncommutative context, wherein it can be adjusted so that the derivative of a product satisfies the Leibniz rule. The relationship with discrete calculus is detailed in Section 7. It is an important background to the content of the paper.  Let $f(x)$ denote a function of a real variable $x,$
and $\tilde{f}(x) = f(x+h)$ for
a fixed difference $h.$ Define the {\em discrete derivative} $Df$ by the formula $Df = (\tilde{f} - f)/h.$ The Leibniz rule is not satisfied.
The formula for the discrete derivative of a product is as shown below: 
\begin{equation}
D(fg) = D(f)g + \tilde{f}D(g).
\end{equation}
We can adjust the Leibniz rule by introducing a noncommutative operator $J$ with the property that 
\begin{equation}
fJ = J\tilde{f}. 
\end{equation}

Define a modified discrete derivative by the formula
\begin{equation}
\nabla(f) = JD(f).
\end{equation}

 It follows at once that 
 \begin{equation}
\nabla(fg) = JD(f)g + J\tilde{f}D(g) = JD(f)g + fJD(g) = \nabla(f)g + f\nabla(g).
\end{equation}
Note that 
\begin{equation}
\nabla(f) = (J\tilde{f} - Jf)/h = (fJ-Jf)/h = [f, J/h].
\end{equation}
The modified discrete derivatives are represented by commutators, and satisfy the Leibniz rule. 
Discrete calculus can be embedded into a noncommutative
calculus based on commutators. With this understanding of the relationship of discrete calculus and commutator calculus it is possible to consider discrete models for the structures described in this paper, and it is possible to compare the commutators that arise from discrete observations with the commutators in quantum mechanics.\\

Let  $[A,B] = AB -BA$ denote a commutator in an abstract algebra. Define 
$DA = [A,J]$ for a given element $J.$ Then $D$ is a derivation in the sense that
$D(AB) = D(A)B + AD(B)$ (the Leibniz rule). Once we have derivations,
geometric concepts become available. If two
derivations $\nabla_{J}A = [A,J]$ and $\nabla_{K}A =[A,K]$ are given, then we can form their commutator
\begin{equation}
 [\nabla_{J}, \nabla_{K}]A = \nabla_{J}\nabla_{K}A - \nabla_{K}\nabla_{J}A = [[J,K],A].
\end{equation}
(The verification of this last inequality is given in the next Section.)
$R_{JK} = [J,K]$ is defined to be the {\it curvature} associated with $\nabla_J$ and $\nabla_K.$ The commutator of the derivations $\nabla_J$ and $\nabla_K$ 
is represented by $R_{JK}=[J,K].$ When $J$ and $K$ commute, then the 
derivations themselves commute and the curvature vanishes.
We shall demonstrate that curvature in this sense is the formal analog of the curvature of a gauge connection.
\bigbreak

This paper consists in 9 sections including the introduction.  Section 3 outlines the general properties of 
calculus in a noncommutative domain ${\cal N}$, where derivatives are represented by commutators.  Included is a special element $H$ such that the total time derivative is given by the formula  $\dot{F} = [F,H]$ for any $F$ in the noncommutative domain. We assume an initial ``flat" coordinatization where algebra elements $X^{1},\cdots,X^{n}$ represent position coordinates  and commute with each other, and another set of elements $P_1 , \cdots , P_n $ that commute with one another represent their partial derivatives.  We take  $[X^i , P_j ] = \delta_{ij}$ where this means that the commutator is equal to $0$ unless
$i=j$ when it is equal to $1.$ Hence we can define $\partial F/\partial X^i = [F, P_i ].$ \\
 
A formal analog of Hamilton's equations arises in a 
flat coordinate system, and the Heisenberg version of Schr\"{o}dinger's equation arises as well.  Section 4  explores the consequences of 
defining dynamics in the form 
\begin{equation}
m\dot{X^i} = m dX^{i}/dt = {\cal G}_{i}.
\end{equation}
Here $m$ is a {\it constant}, meaning that $m$ is in the center of the algebra so that it commutes with all elements of the algebra. Since $m$ is the analog of mass, we assume that $m$ is invertible and non-zero. We take $\{ {\cal G}_{1},\cdots, {\cal G}_{d} \}$
to be a collection of elements of the noncommutative domain $\cal N.$ Let ${\cal G}_{i} = P_{i} -  A_{i}.$
This is a definition of $A_{i}$ with $\partial_i F= \partial F/\partial X^{i} = [F,P_{i}].$ Gauge theory formalism appears
via the curvature of $\nabla_i$ with  $\nabla_{i}(F) = [F, P_i - A_i].$  
\begin{equation}
R_{ij} = \partial_{i} A_{j} - \partial_{j} A_{i} + [A_{i}, A_{j}].
\end{equation}
With $m\dot{X^{i}} = {\cal G}_{i}= P_{i} -  A_{i},$ define $ g^{ij} = [X^{i}, m\dot{X^{j}}].$  This is a natural choice for
a generalized metric. For a quadratic Hamiltonian with these metric coefficients, one can prove the formula  $ g^{ij} = [X^{i}, m\dot{X^{j}}].$  
Onc can show that for any $F,$ 
\begin{equation}
\dot{F} = \frac{1}{2} (\dot{X^i} \partial_i (F) + \partial_{i} (F) \dot{X^i}).
\end{equation}
Note that except for the appearance of two orderings of a product, this is
the standard formula for the total derivative in multi-variable calculus. It can happen that, under constraints,  certain basic formulas go directly over to corresponding (ordered) formulas in the noncommutative world.
In this case we see that a quadratic Hamiltonian has this property.
See also  \cite{Constraints}.
A covariant version of the Levi-Civita connection is a consequence. This connection satisfies the formula
\begin{equation}
\Gamma_{kij} + \Gamma_{ikj} = \nabla_{j}g^{ik} = \partial_{j} g^{ik} + [g^{ik}, A_j].
\end{equation}
and so corresponds, in the noncommutative world, to the connection of Hermann Weyl in his original gauge theory \cite{Weyl}. See Section 4 and Section 5 of the present paper.\\
  
Section 6 is a discussion  on the structure of the Einstein tensor and how the Bianchi identity can be seen from the Jacobi idenity in a noncommutative world.  Section 7 is a discussion about how discrete calculus embeds in noncommutative calculus. This section can be regarded as an indication of an arena of applications of the methods of the present paper. In particular, we discuss a model for discrete measurement and show how commutators arise in this model, and how a commutator of position and momentum is solved in the discrete context by a Brownian walk. Section 8 is an exposition on the Weyl $1$-form leading to electromagnetism, its relationship with the Feynman-Dyson derivation of electromagnetismm  and a reminder of how this formalism is generalized to gauge theories, loop quantum gravity and low-dimensional topology. We consider a question about the Ashtekar variables, loop quantum gravity and their relationship with noncommutative worlds.
This question will be taken up in a sequel to the present paper. This section ends with a recaptulation of our derivation \cite{KN:QEM,QG,NCW} of a generalization of the Feynman-Dyson  \cite{Dyson} derivation of electromagnetism from commutator calculus. Our generalization can be compared with Weyl's orginal derivation using differential forms and we plan subsequent work on this aspect.\\

A significant structural point comes forth in Section 8 where we review our generalization of the Feynman-Dyson work. In our approach to this we base the whole derivation on writing deriviatives 
as commutators and demanding a noncommutative world analog of the formula
\begin{equation}
\dot{F} = \partial_{t}F + \Sigma_{i}\dot{X^{i}}\partial_{i}(F)
\end{equation}
See equation (203) for the specifics. The generalized version of electromagnetism follows entirely from this constraint. See also \cite{Constraints,Deakin}. Constraints between the form of the calculus in standard worlds and the form of the calculus in noncommutative worlds seem to be at the heart of physical laws. This needs better understanding.\\

We have taken the liberty of writing the last four sections of this paper to show background ideas and structures that are related to the main themes of the paper and to indicate further lines of inquiry. The Section 6 on the Bianchi identity shows how another aspect of differential geometry appears in the context of commutators.  
The Section 7 on discrete calculus shows the beginning of how the methods of this paper can be applied to discrete formulations of physics. Section 8 on the Weyl
$1$-form gives background for understanding both our noncommutative calculus and the geometric form of electromagnetism. We describe how the Weyl 1-form gives rise to electromagnetism 
derived from a vector and scalar potential. We recall how Weyl's form was generalized to gauge theory where the field is given by the gauge curvature and that this is the local holonomy of the 
generalized Weyl 1-form. We then recall our Electromagnetic Theorem \cite{QG} where we derive a generalization of the Feyman-Dyson work that fits a gauge theory. The analogs of the electric and magnetic fields can then be compared with their counterparts from the Weyl 1-form. We find a startling match that leads to new research problems, as the reader will see on examining the end
of Section 8.\\

Section 9 is a summary of ideas and results, a discussion of further work and references to current work that we feel is relevant to this research.\\

\noindent {\bf Remark.} Our papers \cite{Constraints,Deakin} discuss higher order constraints. Our paper  \cite{KN:QEM} was inspired by the Feynman-Dyson derivation of electromagnetism from commutator calculus \cite{Dyson,Hughes,Mont,Tanimura}.  Other relevant papers are \cite{KN:Dirac,Twist,NonCom,ST,Aspects,QG,Boundaries,NCW,Constraints,Deakin}.

\section{Noncommutative Geometry}
This section is a very concise introduction to the Noncommutative Geometry of Alain Connes  \cite{Connes}. Connes work is relevant to the ideas in this paper. It will be a separate project
to make the comparisons in detail.\\

The classical case is that of commutative ${\cal C}^{*}$ algebras ${\cal A} = {\cal C}(M)$, the algebra of continuous complex valued functions on a compact manifold $M.$ In this case one has,
for physics that $M=T^{*}N$  where $N$ is the cotangent bundle of the space of configurations with its canonical symplectic structure. Here $M$ is the phase space of the physical system.
The Gelfand-Naimark Theorem \cite{Connes} shows that there is an equivalence between the geometric physics of this phase space and  the algebraic persepecitive working only with the space of functions that constitutes the ${\cal C}^{*}$ algebra ${\cal A} .$ Connes takes this Theorem as the lead for studying noncommutative spaces and/or ${\cal C}^{*}$  algebras, by 
defining the relevant calculus and analysis directly in terms of the noncommutative spaces.\\

Here is a direct quote from {\it Chapter IV - Quantized Calculus} in the book by Connes \cite{Connes}.: \\

\noindent "The basic idea of this chapter, and of noncommutative differential geometry
is to {\it quantize} the differential calculus using the following operator theoretic notion for the differential 
\begin{equation}
df = [ F,f].
\end{equation}
Here $f$ is an element of an involutive algebra ${\cal A}$ of operators in a Hilbert space ${\cal H}$, while $F$ is a selfadjoint operator of square one ($F^2 = 1$) in ${\cal H}.$ At first
one should think of $f$ as a function on a manifold, i.e. of ${\cal A}$  as an algebra of functions, but one virtue of our construction is that it will apply in the noncommutative case as well."\\

``Since the word {\it quantization} is often overused we feel the need to justify its use in our context."\\

``{\it First}, in the case of manifolds the above formula replaces the differential $df$ by an operator theoretic expression involving a {\it commutator}, which is similar to the replacement of the
Poisson brackets of classical mechanics by commutators."\\

``{\it Second}, the {\it integrality} aspect of quantization (such as the integrality of the energy levels of the harmonic oscillators) will have as a counterpart the integrality of the index of a 
Fredholm operator, which will play a crucial role in our context."\\

Connes works in general with a noncommutative ${\cal C}^{*}$  algebra ${\cal A}$ of operators on a Hilbert space ${\cal H}.$ The symbol $df$ for him must satisfy
$d^2 f = 0.$ This can be accomplished by assuming that $F$ is in the center of the algebra ${\cal A}$ and the differential $df = [ F,f].$ is a graded differential. This means that 
successive applications of $d$ alternate between commutator brackets $[a,b] = ab-ba$  and mutator brackets $\{a,b\} = ab+ba.$ Then we have
\begin{equation}
d^2 f = \{F, [F, f]\} = F(Ff - f F) + (Ff - fF)F
\end{equation}
\begin{equation}
= FFf - FfF + FfF -fFF = [F^2, f].
\end{equation}
Thus if $F^2$ belongs to the center of the algebra ${\cal A},$ then $d^2 f = 0$ for all $f.$\\

One can refer to $df$ as an {\it infinitesimal}. An operator $T$ on the Hilbert space ${\cal H}$ is said to be infinitesimal if it is {\it compact} where compactness of the operator
is a condition on the eigenvalues of the associated operator $(TT^{*})^{1/2}.$ We will refer to Connes \cite{Connes} for the details of the definition of compactness. 
Differential calculus in the Connes NCG is then given by a triple $({\cal A}, {\cal H}, F)$ where $[F, f]$ is a compact operator for every $f$ in ${\cal A}.$ Such triples are 
called {\it Fredholm Modules.} With this definition, Connes can define higher differentials, Grassmann calculus, cohomology in terms of differential forms and in general lift
quantum physical structure into the category of the Fredholm modules.\\

Connes defines a quantized calculus via {\it derivations} $D: {\cal A} \longrightarrow {\cal E}$ where ${\cal E}$ is an ${\cal A}$ bimodule, and a derivation is a map of modules that satisfies
the Leibniz Rule:
\begin{equation}
D(ab) = D(a)b + aD(b).
\end{equation}
Commutators with fixed elements are examples of such derivations and are called in this theory {\it inner derivations}. In this context one can define universal $n$-forms and the 
appropriate differential graded algebra $\Omega {\cal A} = \bigoplus_{n \in N} \Omega^{n} {\cal A}.$\\

 \noindent {\bf Remarks.} We have included this very skeletal description of the Connes framework of noncommutative geometry to indicate an important context in which it is possible to
 perform physics based on underlying noncommutativity. Note that by making the constructions relative to a ${\cal C}^{*}$ algebra ${\cal A}$ the Connes theory has available selfadjoint operators
 for observables in quantum mechanics, appropriate measure theory and the use of algebraic formulations of calculus. The fact that such constructions are possible can help orient the reader of the present paper, where we concentrate only on structure related to algebraic calculus. Our aims in this paper are different from the aims of the full noncommutative geometry. We examine the
relationship between the standard worlds of smooth calculus and the algebraic worlds with (using the terminology above) inner derivations. We are particularly interested in how simple relationships between the rules of the calculus in the standard worlds and the corresponding rules for the calculus in the noncommutataive worlds affect the analogs of physical equations.
This investigation begins in the next section of the paper. We expect that some of our observations will be of use in the context of noncommutative geometry.\\

\section {Calculus in Noncommutative Worlds}

Let $\cal N$ be an abstract associative algebra  that admits commutators. If $A$ and $B$ are
in $\cal N,$ then $[A,B] = AB - BA$ is also an element of  $\cal N.$
\bigbreak

For a fixed $N$ in $\cal N$ define
\begin{equation}
\nabla_N : \cal N \longrightarrow \cal N 
\end{equation}
by the formula
\begin{equation}
\nabla_{N} F = [F, N] = FN - NF.
\end{equation}
$\nabla_N$ is a derivation satisfying the Leibniz rule. 
\begin{equation} 
\nabla_{N}(FG) = \nabla_{N}(F)G + F\nabla_{N}(G).
\end{equation}
\bigbreak

\noindent Such derivations do not, in general commute with one another. The key result for their non-commutation is as follows.\\

\noindent {\bf Theorem 2.1.} With the definitions as above, the commutator of the two derivations is given by the formula 
\begin{equation}
[\nabla_ J, \nabla_ K]F= [[J,K], F].
\end{equation}
\noindent {\bf Proof.}
\begin{equation}
[\nabla_ J, \nabla_ K]F
\end{equation}
\begin{equation}
(\nabla_ J \nabla_K - \nabla_ K \nabla_J)F 
\end{equation}
\begin{equation}
= [[F,K],J] - [[F,J],K] 
\end{equation}
\begin{equation}
= (FK-KF)J - J(FK-KJ) - (FJ-JF)K + K(FJ-JF)
\end{equation}
\begin{equation}
= FKJ-KFJ-JFK+JKF-FJK+JFK+KFJ-KJF
\end{equation}
\begin{equation}
= (JK-KJ)F - F(JK - KJ)
\end{equation}
\begin{equation}
= [[J,K], F]
\end{equation}
This completes the proof. $\hfill \Box$\\

$R_{JK} = [J,K]$ is defined to be the {\it curvature} associated with $\nabla_J$ and $\nabla_K.$\\

Within $\cal N$ we build 
a world that imitates the
behaviour of flat coordinates in Euclidean space.  We need that the derivations for those coordinate directions commute with one another.
Suppose that $X$ and $Y$ are coordinates and that $P_X$ and $P_Y$ represent derivatives in these directions so that one writes 
\begin{equation}
\partial_X F = [F, P_X] 
\end{equation}
and   
\begin{equation}
\partial_Y F = [F, P_Y].  
\end{equation}
 
In general, the two derivatives will not commute, but we have that
\begin{equation}
 \partial_ X \partial_Y = \partial_ Y \partial_X
 \end{equation}
 when $[P_X, P_Y] = 0.$   
\bigbreak

\noindent 
Let $X^1, X^2, \cdots, X^d$  represent coordinates.  The $X^i$ satisfy the commutator equations below with the $P_j$ chosen to represent differentiation with
respect to $X^j.$:
\begin{equation}
[X^{i}, X^{j}] = 0
\end{equation}
\begin{equation}
[P_{i},P_{j}]=0
\end{equation}
\begin{equation}
[X^{i},P_{j}] = \delta_{ij}.
\end{equation}
Derivatives are represented by commutators. 
\begin{equation}
\partial_{i}F = \partial F/\partial X^{i} = [F, P_{i}],
\end{equation}
\begin{equation}
\hat{\partial_{i}}F = \partial F/\partial P_{i} = [X^{i},F]. 
\end{equation}

\noindent {\bf Remark on Poisson Brackets.} Consider the case of a single $X$ and a single $P.$ We then have 
$XP-PX=1$ and 
\begin{equation} \partial F/\partial X = [F,P] \end{equation}
and
\begin{equation} \partial F/\partial P = [X,F]. \end{equation}
This is in exact analogy with the Poisson bracket in standard variables $x$ and $p$ where we have
\begin{equation}  \{f,g\} = (\partial f/\partial x)  (\partial g/\partial p) -  (\partial g/\partial x)  (\partial f/\partial p) \end{equation}
and we certainly have that
\begin{equation} \partial f/\partial x = \{x,p\} \end{equation}
and
\begin{equation} \partial f/\partial p = \{x,f\}. \end{equation}
By choosing flat local coordinates as we have done above, we insure that there is a direct correspondence between the noncommutative 
calculus and the standard calculus at this point of departure from Poisson brackets. As soon as we introduce other noncommutative elements into
the algebraic world, the two calculi will diverge from one another. In the discussion below, we will point out that the noncommutative version of a quadratic Hamiltonian
keeps the calculi in close correspondence.\\

\noindent {\bf Introducing Time.}
The time derivative is represented by commutation with $H$ as shown below:

\begin{equation}
dF/dt = [F, H].
\end{equation}
$H$ corresponds to the classical Hamiltonian  or to the Hamiltonian operator in quantum physics. In the abstract world ${\cal N}$ it is neither of these. We can consider representations of 
the algebra ${\cal N}$ where connections with the physics are more direct. For example, if ${\cal N}$ is represented to a ${\cal C}^{*}$  algebra, then the Hamiltonian can be represented
as a self-adjoint operator for the quantum mechanical version where we take $i\hbar dA/dt = [A, H].$\\

{\it Hamilton's equations are a consequence of these definitions.} 
\bigbreak

\noindent {\bf Hamilton's Equations.} 
\begin{equation}
dP_{i}/dt  = -\partial H/\partial X^{i}
\end{equation}
\begin{equation}
dX^{i}/dt = \partial H/\partial P_{i}. 
\end{equation}

\noindent {\bf Proof.}

\begin{equation}
dP_{i}/dt = [P_{i}, H] = -[H, P_{i}] = -\partial H/\partial X^{i}
\end{equation}
\begin{equation}
dX^{i}/dt = [X^{i}, H] = \partial H/\partial P_{i}. 
\end{equation}

This completes the proof.  $\hfill \Box$\\

\noindent {\bf Remark.} These are exactly Hamilton's equations of motion. The pattern of
Hamilton's equations is built into the system. It is natural to ask how this formal appearance of Hamilton's equations is related to their role in physics.
Recall the one variable case. In the standard world we have a coordinate $x$ and the momemtum is given by the formula $p = m dx/dt$ and the energy of the system is
$H = (1/2m)p^2 + V(x).$  Newton's law asserts that $m \ddot{x} = F = - \partial V/\partial x = - \partial H/\partial x.$  Now we have
$dp/dt = m \ddot{x}.$ Thus 
\begin{equation}
dp/dt = - \partial H/\partial x  
\end{equation} 
and
\begin{equation}
dx/dt = p/m = \partial H/\partial p.
\end{equation}
Thus we have Hamilton's equations as an expression of the form of the Hamiltonian in the presence of Newton's law.
Hamilton went on to observe that if $F = F(x,p)$ is a function of position and momentum, then 
\begin{equation}
\dot{F} = (dx/dt) (\partial F/\partial x) + (dp/dt) (\partial F/ \partial p)
\end{equation}
and hence
\begin{equation}
\dot{F} = (\partial F/\partial x)  (\partial H/\partial p) -( \partial H/\partial x)  (\partial F/ \partial p) = \{ F, H \}
\end{equation}
Where the last expression is by definition the {\it Poisson bracket} of $F$ with respect to $x$ and $p.$
The dynamical equation 
$$\dot{F} = \{ F, H \}$$
is the beginning of Hamiltonian Mechanics, and it turns out that the Poisson bracket has the same formal properties as a commutator. We will not go further into Hamiltonian mechanics, but
it should be clear to the reader that what has happened in our noncommutative world formulation is that we developed Hamilton's idea from the other end, starting from a commutator 
analog of Hamilton's dynamical equation, and the Hamilton equations come out as a consequence. It is striking that they do happen when we begin with the concept of derivatives
as commutators.
\bigbreak

\noindent {\bf Remark on Curvature and Covariant Derivative.} By defining $R_{JK} = [J,K]$ as a {\it curvature} associated with the derivations $\nabla_J$ and $\nabla_K$  we are establishing terminology for use
in the algebraic context of noncommutative worlds. Since we have now a reference flat world with derivatives $\partial_j F= \partial F /\partial X^{j}  = [F, P_j],$ we can also create a terminology for (generalized)  {\it covariant derivatives} by defining 
\begin{equation}
\nabla_{j} F = [F, P_j - \Gamma_j] = \partial_j F - [F, \Gamma_j]
\end{equation}
where the $\Gamma_i$ are some given elements of the algebra ${\cal N}.$ Then the $\Gamma_i $ are analogous to a Christoffel symbols and the $\nabla_{j}$ represented by $P_j - \Gamma_j$ can have their curvatures (in the sense of their commutators) analyzed for the purposes of a problem at hand.  \\

The reader may wish to recall that in standard differential geometry (in local coordinates) the curvature  tensor  $R$ is obtained from the commutator of a covariant derivative $\nabla_k,$
associated with the a connection (expressed here in the formalism of a Christoffel symbol)  $\Gamma^{i}_{jk} = (\Gamma_{k})^{i}_{j}.$  
One has , for a vector field $\lambda$ (with local vector coordinates $\lambda_{i}$) the expression for the covariant derivative $\nabla_{j}$, given by a combination of standard differentiation and the application of the connection as a linear transformation.
\begin{equation}
\nabla_{j}\lambda_{i} = \partial_{j}\lambda_{i} - \Gamma^{k}_{ij}\lambda_{k}
\end{equation}
 or
\begin{equation}
\nabla_{j}\lambda = \partial_{j}\lambda - \Gamma_{j}\lambda
 \end{equation}
for a vector field $\lambda.$ 
With 
\begin{equation}
R_{ij} = [\nabla_{i}, \nabla_{j}] = \partial_{j}\Gamma_{i} - \partial_{i}\Gamma_{j} + [\Gamma_{i}, \Gamma_{j}], 
\end{equation}
Here the commutator $[\Gamma_{i}, \Gamma_{j}]$ is the commuator of the linear transformations. The analogy with the noncommutative world is very strong, but 
{\it we are not giving a strict correspondence between curvatures as commutators in the noncommutative world and curvatures as commutators in standard differential geometry.}
The point of our work is to explore the noncommutative world by examining analogous structures that exist within it. Once one takes a flat world for reference, one has the notion of covariant 
derivatives, in our generalized sense, as derivations represented by elements of the form $P_{j} - \Gamma_{j}$ where $\Gamma_{j}$ is any element of the algebra ${\cal N}.$ 
 \\

\section {Dynamics, Gauge Theory and the Weyl Theory.} 
One can take the general dynamical
equation in the form 
\begin{equation}
m \dot{X^i} = m dX^{i}/dt = {\cal G}_{i}
\end{equation}
 where $m$ is a constant (meaning that $m$ is a non-zero, invertible element of the center of the algebra ${\cal N}$) and $\{ {\cal G}_{1},\cdots, {\cal G}_{d} \}$
is a collection of elements of $\cal N.$ Write ${\cal G}_{i}$
relative to the flat coordinates via ${\cal G}_{i} = P_{i} -  A_{i}.$
This defines $A_{i}$ with $\partial_i F= \partial F/\partial X^{i} = [F,P_{i}].$  We define derivations corresponding to the ${\cal G}_{i}$ by the formulas  
\begin{equation}
\nabla_{i}(F) = [F, {\cal G}_{i}],
\end{equation}
then one has
the curvatures (as commutators of derivations) 
\begin{equation}
[\nabla_{i}, \nabla_{j}]F =[ [{\cal G}_{i}, {\cal G}_{j}],F]
\end{equation}
 (by the formula in the introduction to this paper)
 \begin{equation}
= [ [P_{i} -  A_{i}, P_{j} -  A_{j}],F]
\end{equation}
\begin{equation}
 = [[P_i ,P_j]  - [P_i, A_j] - [A_i, P_j] + [A_i, A_j],F]
\end{equation}
\begin{equation}
= [\partial_i A_j - \partial_j A_i + [A_i, A_j], F].
\end{equation}
Thus the curvature is given by the formula 
\begin{equation}
R_{ij} = \partial_{i} A_{j} - \partial_{j} A_{i} + [A_{i}, A_{j}].
\end{equation}
This curvature formula is the noncommutative world analog of the curvature of a gauge
connection when $A_i$ is interpreted as such a connection.  
\bigbreak

With $m \dot{X^{i}} = P_{i} -  A_{i},$ the commutator $[X^{i}, m \dot{X^{j}}]$ takes the form
\begin{equation}
g^{ij}=[X^{i}, m \dot{X^{j}}] = [X^{i}, P_{j} -  A_{j}] = [X^{i}, P_{j}] - [X^{i}, A_{j}] = \delta_{ij} - [X^{i}, A_{j}] .
\end{equation}
Thus we see that the ``gauge field" $A_{j}$ provides the deviation from the Kronecker delta in this commutator. 
We have $[m\dot{X^{i}}, m\dot{X^{j}}] = R_{ij},$ so that these commutators represent the curvature.\\

Before proceeding further it is necessary to explain why $g^{ij}=[X^{i}, m\dot{X^{j}}]$ is a correct analogue to the metric in the classical physical situation.
In the classical picture we have $ds^2 = g^{ij} dX^{i} dX^{j}$ representing the metric. Hence 
\begin{equation}
(m/2)ds^2 /dt^2 = (1/2m)g^{ij} (mdX^{i}/dt)( mdX^{j}/dt ) = (1/2m)g^{ij} p_i p_j 
\end{equation}
and so the 
Hamiltonian is $H = \frac{1}{2m} g^{ij} p_i p_j .$   By convention we sum over repeated indices.\\

We choose an $H$ in the algebra ${\cal N}$ to represent
the total time derivative so that $\dot{F} = [F,H]$ for any $F.$  Let there be given elements $g^{ij}$ such that
\begin{equation}
[g^{ij}, X^{k}]=0
\end{equation}
and 
\begin{equation}
g^{ij} = g^{ji}. 
\end{equation}
{\it Note that the $g^{ij}$ are then formally functions of the space variables $X^{i}$ and do not involve the $P_{j}.$}
This is a natural assumption in analogy to a classical metric where it is a function of the position coordinates.\\

We choose
\begin{equation}
H = \frac{1}{2m} g^{ij}P_{i}P_{j}.  
\end{equation}
This is the noncommutative analog of the classical
$H = \frac{1}{2m} g^{ij}p_{i}p_{j}.$  We now show that this choice of Hamiltonian implies that 
\begin{equation}
[X^{i}, m\dot{X^{j}}] = g^{ij}.
\end{equation}
\\

\noindent {\bf Quadratic Example.} There is an advantage to choosing a quadratic Hamiltonian in the noncommutative world. The next two lemmas show that with this choice, the noncommutative
calculus interacts harmoniously with the patterns of the standard calculus. To show how this works, here is a minature example of what we are about to do in general. Take a noncommutative world with two generating elements $X$ and $P$ and assume that 
\begin{equation} 
[X,P] = XP - PX = 1.
\end{equation}
Let 
\begin{equation} 
H = P^2/2
\end{equation}
and define
\begin{equation} 
\dot{F} = [F,H] = [F, P^2/2].
\end{equation}
Then
\begin{equation} 
\dot{X} = [X,P^2/2] = (1/2)((XP)P - P(PX))
\end{equation}
\begin{equation} 
= (1/2)((1+ PX)P - P(XP -1)) = P.
\end{equation}
Thus
\begin{equation} 
\dot{F} = (1/2)[F,P^2] = (1/2)(FPP - PPF)
\end{equation}
\begin{equation} 
= (1/2)((FP-PF)P + PFP - PPF)
\end{equation}
\begin{equation} 
= (1/2)([F,P] P + P[F,P])
\end{equation}
\begin{equation} 
= (1/2)((\partial F/ \partial X) \dot{X} + \dot{X} (\partial F / \partial X)).
\end{equation}
Thus in this small example we see that the choice of quadratic Hamiltonian leads to the general validity of a symmetrized version of the usual chain rule for 
differentiation with respect to time. We now give a general version of this derivation. By examining the small example, the reader can see that this harmony between
standard calculus and the noncommutative calculus needs a quadratic Hamiltonian.\\

\noindent {\bf Lemma 3.1.} Let $g^{ij}$ be given such that $[g^{ij}, X^{k}]=0$ and $g^{ij} = g^{ji}.$ Define 
\begin{equation}
H = \frac{1}{2m} g^{ij}P_{i}P_{j}
\end{equation}
 (where we sum over the repeated indices) and 
 \begin{equation}
 \dot{F} = [F,H]. 
 \end{equation}
 Then 
\begin{equation}
[X^{i}, m\dot{X^{j}}] = g^{ij}.
\end{equation}
\bigbreak

\noindent {\bf Proof.} Note that $m$ is a non-zero element in the center of the algebra ${\cal N}$, and so may be moved freely in and out of the commutators in the 
calculations performed below. Since
\begin{equation}
\dot{X^k}  = [ X^k , H],
\end{equation}
we have
\begin{equation}
2m \dot{X^k}  = [X^{k}, g^{ij}P_{i}P_{j}] = g^{ij}[X^{k}, P_{i}P_{j}]
\end{equation}
\begin{equation}
= g^{ij}([X^{k}, P_{i}]P_{j} + P_{i}[X^{k}, P_{j}])
\end{equation}
\begin{equation}
= g^{ij}(\delta_{ki} P_{j} + P_{i} \delta_{kj}) = g^{kj}P_{j} + g^{ik}P_{i}
\end{equation}
\begin{equation}
= 2 g^{kj}P_{j}. 
\end{equation}
Thus we have shown that 
\begin{equation}
m \dot{X^k} = g^{kj}P_{j}.
\end{equation}
Then
\begin{equation}
[X^{r}, m\dot{X^{k}}] = [X^{r}, g^{kj}P_{j}] = g^{kj}[X^{r}, P_{j}] = g^{kj} \delta_{rj} = g^{kr} = g^{rk}.
\end{equation}
This completes the proof. $\hfill \Box$
\bigbreak

\noindent {\bf Remark.}  It is worth noting that if we had defined $g^{ij}$ by the formula, $g^{ij} = [X^{i}, m\dot{X^{j}}]$ and assumed that $[g^{ij}, X^{k}]=0,$  then it follows from 
differentiating $[g^{ij}, X^{k}]$ with respect to time that $g^{ij} = g^{ji}.$\\

\noindent {\bf Remark.} We can generalize the form of this Hamiltonian to
\begin{equation}
 H = \frac{1}{4m}(g^{ij}P_{i}P_{j} + P_{i}P_{j}g^{ij}) + V
 \end{equation}
  where $V$ commutes with the $X^{i}.$\\ 

We can further remark that \\

\noindent {\bf Lemma 3.2.} With the same hypotheses as the previous Lemma and with $F$ any element of the given noncommutative world ${\cal N},$ we have the formula
\begin{equation}
\dot{F} = \frac{1}{2} (\dot{X^i} \partial_i (F) + \partial_{i} (F) \dot{X^i}).
\end{equation}
\\

\noindent {\bf Proof.} Note from the previous Lemma that $m\dot{X^k} =  g^{kj}P_{j}.$
\begin{equation}
\dot{F} = [F,H] = [F, \frac{1}{2m} g^{ij} P_i P_j ] = \frac{1}{2m}g^{ij}[F, P_i P_j]
\end{equation}
\begin{equation}
= \frac{1}{2m}g^{ij} ([F,P_i]P_j + P_i [F,P_j])
\end{equation}
\begin{equation}
= \frac{1}{2m}g^{ij}P_i [F, P_j] + [F, P_i] \frac{1}{2m}g^{ij} P_j
\end{equation}
\begin{equation}
= \frac{1}{2} (\dot{X^i}\partial_i(F) + \partial_{i}(F) \dot{X^i}).
\end{equation}
This completes the proof. $\hfill \Box$\\

Using the quadratic Hamiltonian, we have shown that  the basic time derivative formula in standard worlds
\begin{equation}
\dot{F} = \dot{X^i} \partial_i (F)
\end{equation}
has its correct (symmetrized) noncommutative counterpart. In \cite{Constraints}
we say that, using the quadratic Hamiltonian, the noncommutative world satisfies the {\it first constraint.}\\

Using the commutator $[X^{i}, m\dot{X^{j}}] = g^{ij}$,
one can show \cite{NCW,QG} that  
\begin{equation}
\ddot{X^{r}} = G^{r} + F_{rs}\dot{X^{s}} + \Gamma_{rst}\dot{X^{s}}\dot{X^{t}},
\end{equation}
where $G^{r}$ is a scalar field, $F_{rs}$ is a gauge field and $\Gamma_{rst}$ is the Levi-Civita
connection associated with $g^{ij}$ in the given noncommutative world.
\bigbreak

We now differentiate both sides of
the equation 
\begin{equation}
g^{ij}= [X^{i}, m\dot{X^{j}}], 
\end{equation}
and show how the Levi-Civita connection
appears.\\

\noindent Let $D = d/dt$ in the calculations below.
\bigbreak

The Levi-Civita connection (with covariant derivatives $\nabla_{k}$)

\begin{equation}
\Gamma_{kij} =(1/2)(\nabla_{i}g^{jk}+\nabla_{j}g^{ik}-\nabla_{k}g^{ij})  
\end{equation}

\noindent associated with the $g^{ij}$ comes up almost at once from the
differentiation process described below.  Note that a covariant derivative applied to a metric tensor is
called a {\it non-metricity tensor} in the current literature \cite{Hehl} and is related to metric affine gravity.\\

\noindent To see how this happens, view
the following calculation where 
\begin{equation}
\hat{\partial_{i}}\hat{\partial_{j}}F = [X^{i}, [ X^{j}, F]].
\end{equation}
\bigbreak

\noindent We apply the operator $\hat{\partial_{i}}\hat{\partial_{j}}$ to the second
time derivative of
$X^{k}.$
\bigbreak

\noindent {\bf Lemma 3.3} Let 
\begin{equation}
\Gamma_{kij} =(1/2)(\nabla_{i}g^{jk}+\nabla_{j}g^{ik}-\nabla_{k}g^{ij})
\end{equation}
where 
\begin{equation}
\nabla_{i}(F) = [F, m\dot{X^i}]
\end{equation}
 is the covariant derivative generated by 
$m\dot{X^i} = P_i - A_i .$
Then
\begin{equation}
\Gamma_{kij}=(1/2)\hat{\partial_{i}}\hat{\partial_{j}}m\ddot{X^{k}}.  
 \end{equation}
\bigbreak

\noindent {\bf Proof.} 
\noindent Note that by the Leibniz rule

\begin{equation}
D([A,B]) = \dot{[A,B]} = [\dot{A}, B] + [A, \dot{B}],
\end{equation}
 we have

\begin{equation}
\dot{g^{jk}}= [\dot{X^{j}}, m\dot{X^{k}}] + [X^{j}, m\ddot{X^{k}}].
\end{equation}
Therefore
\begin{equation}
\hat{\partial_{i}}\hat{\partial_{j}}m^2 \ddot{X^{k}} = [mX^{i}, [ X^{j}, m\ddot{X^{k}}]]
\end{equation}
\begin{equation}
= [mX^{i}, \dot{g^{jk}} - [\dot{X^{j}}, m\dot{X^{k}}]]
\end{equation}
\begin{equation}
= [X^{i}, m\dot{g^{jk}}] -  [X^{i}, [m\dot{X^{j}}, m\dot{X^{k}}]]
\end{equation}
(Now use the Jacobi identity $[A,[B,C]] + [C,[A,B]] + [B,[C,A]] = 0.$)

\begin{equation}
= [mX^{i}, \dot{g^{jk}}] + [m\dot{X^{k}}, [X^{i}, m\dot{X^{j}}]]+[m\dot{X^{j}}, [ \dot{X^{k}}, mX^{i}]]
\end{equation}
\begin{equation}
=  - [m\dot{X^i}, g^{jk}] + [m\dot{X^{k}}, [X^{i}, m\dot{X^{j}}]]+[m\dot{X^{j}}, [ \dot{X^{k}}, mX^{i}]]
\end{equation}
\begin{equation}
= \nabla_{i}g^{jk} -\nabla_{k}g^{ij} + \nabla_{j}g^{ik}
\end{equation}
\begin{equation}
 = 2\Gamma_{kij}.
\end{equation}
This completes the proof. $\hfill \Box$
\bigbreak

\noindent {\bf Remark.} This derivation confirms our interpretation of 
\begin{equation}
g^{ij} = [X^i, m\dot{X^j}] = [X^i, P_j] - [X^i, A_j] = \delta_{ij} - \partial A_{j}/\partial P_{i}
\end{equation}
as an abstract form of metric.  Note that there may be no given concept of distance in the noncommutative 
world. This suggests a differential geometry based on non-commutativity and
the Jacobi identity.   At this point, it may be important to compare this formalism with the 
way the geometry works in the Connes theory \cite{Connes}. Certainly in that theory there has been a re-evaluation and reconstruction of differential 
geometry based on a noncommutative calculus, and it would be of great interest to trace the role of the Jacobi identity in the Connes quantized calculus.\\

Note that given 
\begin{equation}
 \Gamma_{kij} =(1/2)(\nabla_{i}g^{jk}+\nabla_{j}g^{ik}-\nabla_{k}g^{ij}),
\end{equation}
  we have
  \begin{equation}
 \Gamma_{ikj} =(1/2)(\nabla_{k}g^{ji}+\nabla_{j}g^{ki}-\nabla_{i}g^{kj}).
\end{equation}
 Hence
 \begin{equation}
\Gamma_{kij} + \Gamma_{ikj} = \nabla_{j}g^{ik} = \partial_{j} g^{ik} + [g^{ik}, A_j].
\end{equation}
The noncommutative world Levi-Civita connection differs from a classical Levi-Civita connection via the use of the covariant derivatives $\nabla_{j} = \partial_{j} + A_j.$  
This formalism can be matched with the Levi-Civita connection in Weyl's theory that combines aspects of general relativity with electromagnetism. 
See \cite{Weyl} Chapter 35, page 290 and \cite{Gauge} p. 85.  In Section 8 we discuss Weyl's original approach to electromagnetism.\\

\noindent {\bf Recalling the Standard Levi-Civita Connection.}
Classical Riemannian geometry begins with the standard Levi-Civita connection. Curvature is defined by parallel displacement.
The infinitesimal parallel
translate of a vector $A$ is given by $A' = A + \delta A$ where
\begin{equation} 
\delta A^{k} = -\Gamma^{k}_{ij}A^{i}dX^{j}
\end{equation}

 \noindent  The Christoffel symbols satisfy the symmetry condition
$\Gamma^{k}_{ij} = \Gamma^{k}_{ji}.$ An inner product is given by the formula
\begin{equation}
<A,B> = g^{ij}A^{i}B^{j} 
\end{equation}
\noindent  To require that this inner product be invariant
under parallel displacement is to require that $\delta(g^{ij}A^{i}A^{j}) = 0.$\\

\begin{equation}
\delta(g^{ij}A^{i}A^{j}) = (\partial_{k}g^{ij})A^{i}A^{j}dX^{k} + g^{ij}\delta(A^{i})A^{j} + g^{ij}A^{i}\delta(A^{j})
\end{equation}
\begin{equation}
= (\partial_{k}g^{ij})A^{i}A^{j}dX^{k} - g^{ij}\Gamma^{i}_{rs}A^{r}dX^{s}A^{j} - g^{ij}A^{i}\Gamma^{j}_{rs}A^{r}dX^{s}
\end{equation}
\begin{equation}
= (\partial_{k}g^{ij})A^{i}A^{j}dX^{k} - g^{ij}\Gamma^{i}_{rs}A^{r}A^{j}dX^{s} - g^{ij}\Gamma^{j}_{rs}A^{i}A^{r}dX^{s}
\end{equation}
\begin{equation}
= (\partial_{k}g^{ij})A^{i}A^{j}dX^{k} - g^{sj}\Gamma^{s}_{ik}A^{i}A^{j}dX^{k} - g^{is}\Gamma^{s}_{jk}A^{i}A^{j}dX^{k}
\end{equation}
\noindent Hence

\begin{equation}
(\partial_{k}g^{ij}) = g^{sj}\Gamma^{s}_{ik} + g^{is}\Gamma^{s}_{jk}.
\end{equation}
\noindent From this it follows that 

\begin{equation}
\Gamma_{ijk} = g^{is}\Gamma^{s}_{jk} = (1/2)(\partial_{k}g^{ij} - \partial_{i}g^{jk} +\partial_{j} g^{ik}).
\end{equation}
 
To generalize the above into the noncommutative context  becomes a significant program for further investigation in the noncommutative world. It would appear that a standard version of this program was implicit in Weyl's original work. See  his papers and the book ``Space Time Matter" \cite{Weyl}.  Our approach suggests a new start on this problem. \\

\section{Recapitulation - Curvature, Jacobi Identity and the Levi-Civita Connection}
We recapitulate and set the stage for a next level of structure.
We use a partially index-free notation. 
Nested subscripts are avoided by using different variable names and then using these 
names in place of  subscripts. We write $X$ and $Y$ instead of
$X^{i}$ and $X^{j}.$ We write $g^{XY}$ instead of $g^{ij}.$ The derivation
$DX$ has the form $DX = [X,J]$ for some $J.$
\bigbreak

$[A,B]$ is assumed to satisfy
the Jacobi identity, bilinearity in each variable, and the Leibniz rule for all functions of the
form $\delta_{K}(A) = [A,K].$

\begin{equation}
\delta_{K}(AB) =  \delta_{K}(A)B + A\delta_{K}(B).
\end{equation}
\bigbreak

\noindent We  consider derivatives in the form 
\begin{equation}
\nabla_{X}(A) = [A, \Lambda_{X}].
\end{equation}
\noindent Examine the following computation:
\begin{equation}
\nabla_{X}\nabla_{Y}F = [[F,\Lambda_{Y}],\Lambda_{X}] = -[[\Lambda_{X}, F],\Lambda_{Y}]- [[\Lambda_{Y},\Lambda_{X}],F]
\end{equation}
\begin{equation}
= [[F, \Lambda_{X}],\Lambda_{Y}] + [[\Lambda_{X},\Lambda_{Y}],F]
\end{equation}
\begin{equation}
= \nabla_{Y}\nabla_{X}F + [[\Lambda_{X},\Lambda_{Y}],F].
\end{equation}
\noindent Thus
\begin{equation}
[\nabla_{X}, \nabla_{Y}]F= R_{XY}F
\end{equation}

\noindent where 
\begin{equation}
R_{XY}F = [[\Lambda_{X},\Lambda_{Y}],F].
\end{equation}
\noindent We can regard $R_{XY}$ as a curvature operator. 
\bigbreak

We assume position variables (operators)
$X$, $Y$, $\cdots$ and momentum variables (operators) $P_X$, $P_Y$, $\cdots$ satisfying the 
equations below.
\begin{equation}
[X,Y]=0
\end{equation}
\begin{equation}
[P_X,P_Y]=0
\end{equation}
\begin{equation}
[X, P_Y] = \delta_{XY}
\end{equation}
\noindent where $\delta_{XY}$ is equal to one if $X$ equals $Y$ and is zero otherwise.
We define 
\begin{equation}
\partial_{X}F = [F,P_X]
\end{equation}
\noindent and 
 \begin{equation}
 \partial_{P_X}F = [X,F]. 
 \end{equation}
\noindent These derivatives
behave correctly in that 
\begin{equation}
\partial_{X}(Y) = \delta_{XY}
\end{equation}
\noindent and 
\begin{equation}
\partial_{P_X}(P_Y) = \delta_{XY} 
\end{equation}
\begin{equation}
\partial_{P_X}(Y) = 0 = \partial_{X}(P_Y)
\end{equation}
 \noindent with the last equations valid 
even if $X=Y.$  \\

With this reference point of (algebraic) flat space we define
\begin{equation} 
\hat{P_{X}} = P_{X} - A_{X}
\end{equation}
\noindent for an arbitrary algebra-valued function of the variable
$X.$  
 \noindent With respect to this deformed momentum 
we have the covariant derivative 
\begin{equation}
\nabla_{X}F = [F,\hat{P_Y}] = [F, P_Y - A_Y] = \partial_{Y}F - [F, A_Y].
\end{equation}
\noindent The curvature for this covariant derivative is given by the formula
\begin{equation}
 R_{XY}F = [\nabla_{X}, \nabla_{Y}]F =  [[\lambda_{X},\lambda_{Y}],F]
 \end{equation}
\noindent where $\lambda_X = P_X - A_X.$  Hence
\begin{equation}
R_{XY} = [P_X - A_X, P_Y - A_Y] = -[P_X,A_Y] - [A_X,P_Y] + [A_X, A_Y]
\end{equation}
\begin{equation}
= \partial_{X}A_Y - \partial_{Y}A_X + [A_X, A_Y].
\end{equation}
\noindent and this has the abstract form of the curvature of a
Yang-Mills gauge field. 
\bigbreak

We compute
\begin{equation}
[X, \hat{P_Y}] = [X, P_Y - A_Y] = \delta_{XY} - [X, A_{Y}].
\end{equation}

\noindent Let 
\begin{equation}
g^{XY} = \delta_{XY} - [X, A_{Y}]
\end{equation}
 \noindent
so that 
\begin{equation}
[X, \hat{P_Y}] = g^{XY}.
\end{equation}
It is useful to restrict to the case where $[X, A_Y]=0$ so that $g^{XY}=\delta_{XY}$ (for the space coordinates).
In order to enter this domain, we set 
\begin{equation}
m\dot{X} = mDX = \hat{P_X} = P_X - A_X. 
\end{equation}
where $m$ is a constant (a non-zero, invertible element of the center of the algebra).
\noindent  We examine the structure of the following special axioms for a bracket.
\begin{equation}
[X, DY] = g^{XY}
\end{equation}
\begin{equation}
[X,Y]=0
\end{equation}
\begin{equation}
[Z,g^{XY}]=0
\end{equation}
\begin{equation}
[g^{XY},g^{ZW}]=0
\end{equation}
Note that 
\begin{equation}
Dg^{YZ} =D[Y,DZ] = [DY,DZ] + [Y, D^{2}Z].
\end{equation}
 \noindent and that $[Z,g^{XY}]=0$ 
implies that
\begin{equation}
[g^{XY},DZ] = [Z, Dg^{XY}]. 
\end{equation}
\noindent
\bigbreak

Define two types of derivations as follows 
\begin{equation}
\nabla_{X}(F) = [F,DX]
\end{equation}
\noindent and 
\begin{equation}
\nabla_{DX}(F) = [X,F].
\end{equation}
 \noindent these are dual with respect to $g^{XY}$ and will
act like partials with respect to these variables in the special case when $g^{XY}$ is a 
Kronecker delta, $\delta_{XY}.$ If the form $g^{XY}$ is invertible, then we can rewrite these
derivations by contracting the inverse of $g$ to obtain standard formal partials.
\bigbreak

\begin{equation}
\nabla_{DX}\nabla_{DY} D^{2}Z = [X,[Y,D^{2}Z]]
\end{equation}
\begin{equation}
= [X,Dg^{YZ} - [DY,DZ]] = [X,Dg^{YZ}] - [X,[DY,DZ]] 
\end{equation}
\begin{equation}
=[g^{YZ},DX] - [X,[DY,DZ]]
\end{equation}
\begin{equation}
=\nabla_{X}(g^{YZ}) -[X,[DY,DZ]].
\end{equation}

\noindent Now use the Jacobi identity on the second term and obtain

\begin{equation}
\nabla_{DX}\nabla_{DY} D^{2}Z = \nabla_{X}(g^{YZ}) + [DZ,[X,DY]] + [DY,[DZ,X]] 
\end{equation}
\begin{equation}
= \nabla_{X}(g^{YZ}) - \nabla_{Z}(g^{XY}) + \nabla_{Y}(g^{XZ}). 
\end{equation}
\noindent This is the formal Levi-Civita connection.
\bigbreak

\section{Einstein's Equations and the Bianchi Identity}
The Bianchi identity (see below for its definition) appears in the context of 
noncommutative worlds as a form of the Jacobi identity. We will explain this, and discuss the classical background \cite{Dirac}. \\

The basic tensor in Einstein's theory of general relativity is 
\begin{equation}
G^{ab} = R^{ab} - \frac{1}{2}Rg^{ab}.
\end{equation}

$R^{ab}$ is the Ricci tensor, $R$ the scalar curvature. These are both obtained
by contraction from the Riemann curvature tensor $R^{a}_{bcd}$ with $R_{ab} = R^{c}_{abc}, R^{ab} = g^{ai}g^{bj}R_{ij},$ and 
$R = g^{ij}R_{ij}.$ Since the Einstein tensor $G^{ab}$ has vanishing divergence, it can be proportional to the
energy momentum tensor $T^{\mu \nu}.$ Einstein's field equations are 
\begin{equation}
R^{\mu \nu} - \frac{1}{2}Rg^{\mu \nu} = \kappa T^{\mu \nu}.
\end{equation}
\bigbreak

The Riemann tensor is obtained from the commutator of a covariant derivative $\nabla_k,$
associated with the Levi-Civita connection $\Gamma^{i}_{jk} = (\Gamma_{k})^{i}_{j}$ (using the space-time metric $g^{ij}$).\\

We can write the formalism in the gauge form by hiding some indices.\\

\begin{equation}
\lambda_{a:b} = \nabla_{b}\lambda_{a} = \partial_{b}\lambda_{a} - \Gamma^{d}_{ab}\lambda_{d}
\end{equation}
 or
 \begin{equation}
\lambda_{:b}= \nabla_{b}\lambda = \partial_{b}\lambda - \Gamma_{b}\lambda
 \end{equation}
for a vector field $\lambda.$ 
With 
\begin{equation}
R_{ij} = [\nabla_{i}, \nabla_{j}] = \partial_{j}\Gamma_{i} - \partial_{i}\Gamma_{j} + [\Gamma_{i}, \Gamma_{j}], 
\end{equation}
one has
\begin{equation}
R^{a}_{bcd} = (R_{cd})^{a}_{b}.
\end{equation}
 (Here $R_{cd}$ is {\it not} the Ricci tensor. It is the Riemann tensor with two internal indices
hidden from sight.)\\

One has explicitly that \cite{Dirac}
\begin{equation}
R_{\mu \nu \rho \sigma }= \frac{1}{2} ( g_{\mu \sigma , \nu \rho} - g_{\nu \sigma , \mu \rho}  - g_{\mu \rho , \nu \sigma} + g_{\nu \rho , \mu \sigma}) + 
\Gamma_{\beta \mu \sigma} \Gamma^{\beta}_{\nu \rho} + \Gamma_{\beta \mu \rho} \Gamma^{\beta}_{\nu \rho}. 
\end{equation}

Symmetries of the Riemann tensor follow from the above formula. When derivatives are replaced by covariant derivatives, symmetries may not
survive.  That is a project for future work.\\

The Bianchi identity states
\begin{equation}
R^{a}_{bcd:e} + R^{a}_{bde:c} + R^{a}_{bec:d} = 0
\end{equation}
 where each index after a colon indicates a covariant derivative.
This can be written in the form 
\begin{equation}
(R_{cd:e})^{a}_{b} + (R_{de:c})^{a}_{b} + (R_{ec:d})^{a}_{b} = 0.
\end{equation}
Bianchi identity follows from local properties of the Levi-Civita connection and symmetries of the Riemann
tensor.A relevant symmetry of Riemann tensor is the equation
$R^{a}_{bcd} = - R^{a}_{bdc}.$

Contraction of the Bianchi identity leads to the Einstein tensor.  
\begin{equation}
R^{a}_{bca:e} + R^{a}_{bae:c} + R^{a}_{bec:a} = 0
\end{equation}
This is he same as
\begin{equation}
R_{bc:e} - R_{be:c} + R^{a}_{bec:a} = 0.
\end{equation}
 Contract once more to obtain
 \begin{equation}
R_{bc:b} - R_{bb:c} + R^{a}_{bbc:a} = 0,
\end{equation}
 and raise indices
 \begin{equation}
R^{b}_{c:b} - R_{:c} + R^{ab}_{bc:a} = 0.
\end{equation}
 Further symmetry gives
 \begin{equation}
R^{ab}_{bc:a} = R^{ba}_{cb:a} = R^{a}_{c:a} = R^{b}_{c:b}.
\end{equation}
 Hence we have
 \begin{equation}
2R^{b}_{c:b} - R_{:c} = 0,
\end{equation}
 which is equivalent to the equation
 \begin{equation}
(R^{b}_{c} - \frac{1}{2}R\delta^{b}_{c})_{:b} = G^{b}_{c:b} = 0.
\end{equation}
 From this we conclude
that $G^{bc}_{:b} = 0.$ 
\bigbreak

\noindent {\bf Bianchi Identity and Jacobi Identity.} Now work in noncommutative worlds. We
have convariant derivatives of the form 
\begin{equation}
F_{:a} = \nabla_{a}F = [F, N_{a}]
\end{equation}
 for elements  $N_{a}$ 
in the noncommutative  world. Choose a covariant derivative. Then we have the curvature
\begin{equation}
R_{ij} = [N_{i},N_{j}].
\end{equation}
Note that $R_{ij}$ is not a Ricci tensor. 
We then have the Jacobi identity
\begin{equation}
[[N_a,N_b],N_c] + [[N_c,N_a],N_b] + [[N_b,N_c],N_a] = 0.
\end{equation}
 Writing Jacobi identity using curvature and covariant
differention we have 
\begin{equation}
R_{ab:c} + R_{ca:b} + R_{bc:a}.
\end{equation}
 In a noncommutative world, every covariant derivative
satisfies its own Bianchi identity.  
\bigbreak

 \section{Discrete Calculus Reformulated with Commutators}

Let $f(x)$ denote a function of a real variable $x,$
Let $\tilde{f}(x) = f(x+h)$ for
some fixed difference $h.$ Define the {\em discrete derivative} $Df$ by the formula $Df = (\tilde{f} - f)/h.$ 
One has the basic formula for the discrete derivative
of a product: 
\begin{equation}
D(fg) = D(f)g + \tilde{f}D(g).
\end{equation}
In discrete calculus the Leibniz rule is not satisfied.
Introduce a new invertible operator $J$ with defining property that  
\begin{equation}
fJ = J\tilde{f}.
\end{equation}
 Define an adjusted discrete derivative by the formula
 \begin{equation}
\nabla(f) = JD(f).
\end{equation}
Then 
 \begin{equation}
\nabla(fg) = JD(f)g + J\tilde{f}D(g) = JD(f)g + fJD(g) = \nabla(f)g + f\nabla(g).
\end{equation}
Note that 
\begin{equation}
\nabla(f) = (J\tilde{f} - Jf)/h = (fJ-Jf)/h = [f, J/h].
\end{equation}
 In the adjusted algebra, discrete derivatives are represented by commutators, and satisfy the Leibniz rule. 
One can see discrete calculus as a subset of a noncommutative
calculus based on commutators. For other relationships with discrete calculus, see \cite{Dimakis}.
\bigbreak

\noindent {\bf Discrete Measurement.} In the noncommutative world, consider a {\it time series} $\{X, X', X'', \cdots \}$ with commuting scalar values.
Let 
\begin{equation}
\dot{X} = \nabla X = JDX = J(X'-X)/\tau
\end{equation}
 where $\tau$ is an elementary time step  The operator $J$ is defined by the equation
$XJ = JX'$ or $X^{t}J = JX^{t + \tau}.$
Moving $J$ across a variable from left to right, is an algebraic model for one tick of the clock.  
\bigbreak

Consider observing $X$ at a given time
and observing (or computing) $DX$ at a given time. 
Since $X$ and $X'$ are parts of computing $(X'-X)/\tau,$ the value associated with $DX,$ the
clock must tick once to find $DX.$ Thus, in measurement, $X$ and $DX$ do not commute.
 
\begin{enumerate}
\item Let $\dot{X}X$ denote the sequence: observe $X$, then obtain $\dot{X}.$ 
\item Let $X\dot{X}$ denote the sequence: obtain $\dot{X}$, then observe $X.$ 
\end{enumerate}
\bigbreak

The commutator $[X, \dot{X}]$ expresses the difference between two orders of discrete measurement.
When the elements of the time series are commuting scalars, one has
\begin{equation}
[X,\dot{X}] = X\dot{X} - \dot{X}X =J(X'-X)^{2}/\tau.
\end{equation}
Thus one can interpret  
\begin{equation}
[X,\dot{X}] = Jk
\end{equation}
 ($k$ a constant ) as 
 \begin{equation}
 (X'-X)^{2}/\tau = k.
 \end{equation}
  The process is a walk with spatial step 
\begin{equation}
\Delta = \pm \sqrt{k\tau}
\end{equation}
 where $k$ is a constant.  
 \begin{equation}
k = \Delta^{2}/\tau.
\end{equation}
Hence $k$  is a diffusion constant for a Brownian walk.
The walk with spatial step  $\Delta$ and time step $\tau$ satisfies the commutator equation above
exactly when  $\Delta^{2}/\tau$ remains constant. The diffusion constant of a Brownian process occurs independent of issues about
probability and continuum limits.  
\bigbreak

\section{On Weyl's 1-form for Electromagnetism and the Feynman-Dyson Derivation of Electromagnetism from Commutators}
In this appendix we review the essentials of Weyl's approach to electromagnetism. This is lucidly explained in \cite{Weyl,Gauge,Eddington,Mori}.
Consider a line element for spacetime of the form 
\begin{equation}
\lambda = F dx + G dy + H dz - \phi dt
\end{equation}
 Regard $\lambda$ as a differential $1$-form.
Then (with wedge products of differentials so that $dx \wedge dy = - dy \wedge dx$ and so on) we have
\begin{equation}
d \lambda = (G_x - F_y) dx \wedge dy + (H_x - F_z) dx \wedge  dz  + (H_y - G_z) dy \wedge dz
\end{equation}
\begin{equation}
 - (F_t + \phi_x)dx \wedge  dt -(G_t + \phi_y)dy \wedge  dt - (H_t + \phi_z) dz \wedge  dt.
\end{equation}
Hence, if we set
\begin{equation}
d \lambda = B_1 dy \wedge dz - B_2 dx \wedge  dz + B_3 dx \wedge dy + E_1 dx \wedge dt + E_ 2 dy \wedge  dt + E_3 dz \wedge dt
\end{equation}
and $\nabla = (\partial/\partial x, \partial/\partial y,\partial/\partial z ),$
${\cal A} = (A_1, A_2, A_3) = (F,G,H)$ then 
\begin{equation}
{\cal E} = - \nabla \phi - \partial {\cal A}/\partial t,
\end{equation}
\begin{equation}
{\cal B} = \nabla \times {\cal A}.
\end{equation}

We refer to Weyl's differential 1-form $\lambda$ as his {\it line element} because, for him it represented a new element in the differential geometry of spacetime.
The 1-form has since found a more coherent place in quantum mechanical contexts.
The differential of the 1-form $\lambda$ produces electric and magnetic fields with the space parts acting as the vector potential and the time part acting as the 
scalar potential. Furthermore, one finds that $d^2 \lambda = 0$ (consequence of the properties of differential forms) and with $d \lambda$ in terms of $\cal{E}$ and $B$ the equation 
$d^2 \lambda = 0$ becomes 
\begin{equation}
\nabla \bullet {\cal B} = 0,
\end{equation}
\begin{equation}
 \nabla \times {\cal E} + \partial {\cal B}/ \partial t = 0.
 \end{equation}
  Thus indeed the line element does represent the potentials for electromagnetism, and the equation $d^2 \lambda = 0$ produces Maxwell's equations. The other two Maxwell equations 
  \begin{equation}
 \nabla \bullet {\cal E} = \rho,
 \end{equation}
 \begin{equation} 
 \nabla {\cal B} - \partial {\cal E} / \partial t = J
\end{equation}
 can be regarded as the definitions of the charge density $\rho$ and the current $J.$\\

This means that we can regard a spacetime line element $\lambda$ as the holder of the structure that gives rise to the electromagnetic field. If $d \lambda = 0$ then the line element will have
no holonomy, no change along different paths from one point to another. But if the ${\cal E}$ and ${\cal B}$ fields defined by $d \lambda$ are non-zero, then distances will vary depending upon the path taken between two points. Thus the curvature of this gauge field was identified by Weyl as the electromagnetic field, and he worked on a formalism to unify it with general relativity. The intuitive idea was that in moving from one point of spacetime to another there was spacetime curvature as in general relativity and also curvature that connoted the electromagnetic field via the variation of the line element. Eventually all these considerations were integrated into physics in a different way by regarding that line element as representing the phase of the quantum wave function. Unifications of gauge theory and general relativity have proceeded in different directions. Here we have begun a different way to formulate the Weyl idea in terms of noncommutative worlds, and it remains to be seen the full consequences of our approach.\\

We remark that the standard generalization of the differential $1$-form $\lambda$ is to write $A = \sum_{i} A_{i} dx^{i}$ as a gauge connection where the $A_{i}$ do not commute with one another and take the form $A_{i}(x) = \sum_{a} A_{i}^{a}(x) T_{a}$ where the $T_{a}$ run over a basis for a matrix representation of the Lie algebra of the gauge group, and the  $A_{i}^{a}(x)$ are smooth
functions on the spacetime manifold. Then the curvature of the gauge connection is $F = dA + A \wedge A$ where $\wedge$ denotes the wedge product of differential forms. This generalizes the 
way we have just described the electromagnetic field in terms of the Weyl differential $1$-form and gives rise to the Yang-Mills fields. At the level of noncommutative worlds we can consider abstract differential forms $A = \sum_{i} A_{i} dx^{i}$ without assuming that the $A_{i}$ are represented in terms of a specific classical gauge group. By the same token, we can examine the structure of covariant derivatives of the form $\nabla_{i} = \partial_{i} + A_{i}$ and indeed one finds directly that 
\begin{equation}
[\nabla_{i}, \nabla_{j}] = \partial_{i}A_{j} - \partial_{j}A_{i} + [A_{i}, A_{j}].
\end{equation}
In this way, the formalism of the differential forms and the formalism of the commutators of covariant derivatives come together.\\

We end this section with a recollection of our previous derivation of a generalization of gauge theory in electromagnetic form via noncommutative worlds.
It is of interest to compare the form of this work with the structure of electormagnetism that comes from Weyl's 1-form.\\

\noindent {\bf Generalizing Feynman-Dyson.}
\begin{enumerate}
\item We do not assume that $[X^{i}, m\dot{X^{j}}] = \delta_{ij},$ nor do we assume $[X^{i},X^{j}]=0.$  
We do assume three coordinate variables $\{ X^1, X^2, X^3\}$ in a given noncommutative world.
\item We define 
\begin{equation}
\partial_{i}(F) = [F, \dot{X^{i}}],
\end{equation}

 and the reader should note
that, these spatial derivations are no longer flat in the sense of our earlier sections (nor were they in the original Feynman-Dyson derivation).
\item Define  $\partial_{t} = \partial/\partial t$ as below
\begin{equation}
 \dot{F} = \partial_{t}F + \Sigma_{i}\dot{X^{i}}\partial_{i}(F) =  \partial_{t}F +  \Sigma_{i} \dot{X^{i}}[F, \dot{X^{i}}].
\end{equation}
 
\item In defining 
\begin{equation}
\partial_{t}F = \dot{F} - \Sigma_{i} \dot{X^{i}}[F, \dot{X^{i}}],
\end{equation}
we use the definition itself to create a
distinction between space and time in the noncommutative world.  
\item The reader can verify the following formula:
\begin{equation}
\partial_{t}(FG) = \partial_{t}(F)G + F\partial_{t}(G) + \Sigma_{i}\partial_{i}(F)\partial_{i}(G).
\end{equation}
$\partial_{t}$ does not satisfy the Leibniz rule in our noncommutative context.
Thus $\partial_{t}$ is an operator that does not have a
representation as a commutator. 
\item Divergence and curl are defined by the equations
\begin{equation}
\nabla \bullet B = \Sigma_{i=1}^{3} \partial_{i}(B_{i})
\end{equation}
and 
\begin{equation}
(\nabla \times E)_{k} = \epsilon_{ijk}\partial_{i}(E_{j}).
\end{equation}
where $\epsilon_{ijk}$ is the well known ``epsilon tensor" that is equal to $+1$ for an even permutation of $123$, $-1$ for an odd permutation of $123$ and $0$ if any two
indices are repeated.  Note that the epsilon tensor obeys the identity 
\begin{equation}
\sum_{i} \epsilon_{abi} \epsilon_{dci} = - \delta^{a}_{c} \delta^{b}_{d} + \delta^{a}_{d} \delta^{b}_{c}.
\end{equation}
\end{enumerate}

\noindent The epsilon identity can be used to rewrite equation (208) as
\begin{equation}
 \dot{F} = \partial_{t}F + \dot{X}  \times ( F \times \dot{X}) + (\dot{X} \bullet F) \dot{X} - (\dot{X} \bullet \dot{X}) F.
\end{equation}

\noindent The last equation follows directly from the work in \cite{QG}. By substituting $\dot{X}$ for $F$ we find the equation 
\begin{equation}
\ddot{X} = \partial_{t}\dot{X} + \dot{X} \times (\dot{X} \times \dot{X}).
\end{equation}
This is our motivation for defining
${\cal E} = \partial_{t}\dot{X}$ and 
${\cal B} = \dot{X} \times \dot{X}.$
With these definitions in place, we have
$\ddot{X} = {\cal E} + \dot{X} \times {\cal B},$ giving an analog of the Lorentz force law for
this theory.  Further calculations yield the following Theorem.\\

\noindent {\bf Electromagnetic Theorem \cite{QG}} With the above definitions of the operators, and taking
\begin{equation}
\nabla^{2} = \partial_{1}^{2} + \partial_{2}^{2} + \partial_{3}^{2}, \,\,\, {\cal B} = \dot{X} \times \dot{X} \,\,\, \mbox{and} \,\,\, 
{\cal E} = \partial_{t}\dot{X} \,\,\, \mbox{we have}
\end{equation}
\begin{enumerate}
\item $\ddot{X} = {\cal E} + \dot{X} \times {\cal B}$
\item $\nabla \bullet {\cal B} = 0$
\item $\partial_{t} {\cal B} + \nabla \times {\cal E} = {\cal B} \times {\cal B}$
\item $\partial_{t}{\cal E} - \nabla \times {\cal B} = (\partial_{t}^{2} - \nabla^{2})\dot{X}$
\end{enumerate}
\bigbreak

\noindent{\bf A Deeper Comparison.}
Now we can go further and compare this Theorem with the Weyl approach via differential forms. For note that $\dot{X} = P - A$ in our formalism where $P$ is the vector 
$P=(P_1,P_2,P_3)$ of representatives for $\partial/\partial X^{i}$ in a flat noncommutative world. In the Electromagnetic Theorem we have 
\begin{equation}
{\cal B} = \dot{X} \times \dot{X}  = (P-A) \times (P-A) = P \times P - (P \times A + A \times P)  + A \times A.
\end{equation}
\noindent Letting
\begin{equation}
 \tilde{\nabla} F = ([F, P_1],[F, P_2],[F, P_3])=(\partial F/\partial X^1, \partial F/\partial X^2,\partial F/\partial X^3),
 \end{equation}
the reader will can check that $P \times P = 0$ and  that 
\begin{equation}
P \times A + A \times P = \tilde{\nabla} \times A.
\end{equation}
And so we have
\begin{equation}
{\cal B} = -\tilde{ \nabla} \times A + A \times A.
\end{equation}
In particular,
if the coefficients of $A$ are commutative (as in standard electtomagnetism), then 
\begin{equation}
{\cal E} = \partial_{t}\dot{X}= \partial_{t}P - \partial_{t}A
\end{equation}
and
\begin{equation}
 {\cal B} = - \nabla \times A.
\end{equation}
This is to be compared with the results from differentiating the Weyl form, the expression of the field in terms of scalar and vector potential:
\begin{equation}
{\cal E} = - \nabla \phi - \partial {\cal A}/\partial t,
\end{equation}
\begin{equation}
{\cal B} = \nabla \times {\cal A}.
\end{equation}

We see that, up to shifiting a sign, the significant point is that $\partial_{t}P$ corresponds to  $\nabla \phi,$ the gradient of the scalar potential. 
In our theory there is no scalar potential, but this correspondence can be explored. We see that the Electromagnetic Theroem is probing the same territory as the Weyl form.
And we see that, since $A \times A$ corresponds to the wedge ${\cal A} \wedge {\cal A},$ this correspondence goes over to the full gauge theory. It will give a new way to understand
the extra appearance of the ${\cal A} \wedge {\cal A}$ in the gauge theory. We usually think of the curvature ${\cal F} = d {\cal A}  + {\cal A} \wedge {\cal A}$ as motivated by the calculation 
of local holonomy of the gauge field. Here it appears inevitably from the structure of non-commutativity. That is in accord with the theme of this paper.\\

Weyl's interpretation of the properties of the line element $A=\lambda$ was that an integral along a path from event $p$ to event $q$, $\int_{p}^{q} A$, would be path dependent and that this 
would represent changes in spacetime distance between points depending on the path (history) between them. This path dependence would be a manifestation of the electromagnetic field
$dA$ (In Weyl's form $A \wedge A = 0$). Einstein criticized the theory on these grounds and a new interpretation eventually appeared. The new interpretation can be summarized by multiplying the integral by the square root of negative unity, $i\int_{p}^{q} A,$ and interpreting it via $e^{i\int_{p}^{q} A}$ as a change of phase of a quantum wave function associated with the gauge connection. More particularly, one puts the electromagnetic potential into the quantum Hamiltonian. Schr{\"o}dinger's equation then  has the form 
\begin{equation}
i \hbar \partial \psi/\partial t = \hat{H} \psi
\end{equation}
where 
\begin{equation}
\hat{H} = -\hat{p}^2/2m + e \phi + V
\end{equation}
 is the Hamiltonian operator, $\hat{p}$ is the canonical momentum operator, $e$ is the electric charge. Here one takes the 
canonical momentum to be given by the formula
\begin{equation}
\hat{p} = -i \hbar \nabla - eA.
\end{equation}

Here we see the reflection of our noncommutative world operator $P - A$ in the standard quantum theory.
In this form, many years later \cite{AB} exactly this effect was discovered for electromagnetism, and it became known as the Aharanov-Bohm effect. The interpretation of the gauge connection 
for phases of quantum wave functions became an established part of physics, vindicating Weyl's intuitions, albeit with a shift of interpretation \cite{Yang}. It is only more recently that gravity is seen in relation to gauge fields. The work of Abhay Ashtekar, Carol Rovelli and Lee Smolin has led to the emergence of the field of Loop Quantum Gravity \cite{Baez,Rovelli,Pullin} where a 
gauge formulation of quantum gravity has non-trivial holonomies for macroscopic loops that are central to the theory. It should be mentioned that in the work of Witten \cite{Witten} these kinds of holonomies are closely related to topological invariants of knots, links and three-manifolds. See also \cite{Kauff:KP}.\\

A particularly interesting theory to examine in our noncommutative context is the loop quantum gravity version of general relativity that uses the Ashtekar variables \cite{Baez,Rovelli,Pullin}. In that theory the metric is expressed in terms of a gauge group and the gauge holonomy plays a significant role in the physics and its relation to topology. We intend to examine this structure in a sequel to the present paper.\\

\noindent {\bf Remark.} In comparing our Electromagnetic Theorem with the Weyl 1-form we see that the simplest, perhaps deepest, mathematical 
commonality is in the presence of the epsilon tensor in both structures. The epsilon appears explicitly in ours via the curl and via the definition of the ${\cal B}$-field.
The epsilon is the same as the fundamental antisymmetry of the Grassmann multiplication of differential forms. We treat time in a special way in our derivatives.
Weyl's 1-form is adjusted to handle a temporal component. This comparison is a beginning for future research.\\
 
\section{Conclusions}
In this paper we have explored calculus based on commutators so that derivations are represented in the form $\nabla_{J}(F) = FJ-JF = [F,J]$ in a given algebra $\cal N$ that is closed under
the operation of commutation. We first noted that $[\nabla_{J}, \nabla_{K}]F = [[J,K], F]$ so that the deviation of our derivations from commutativity is measured by the commutators of the 
operators that represent the derivations. We  defined curvature operators  $R_{JK} = [J,K]$ associated to each such pair of derivations and showed how the formalism of the non-commuttive 
calculus aligns itself with physics. A flat framework for physics can be constructed by taking a collection of position coordinates $\{X^{i}\}$ that commute with each other and a collection
of operators $\{P_{i}\}$ that also commute with one another, and we assume that $[X^{i}, P_{j}] = \delta_{ij}$ where $\delta_{ij}$ is equal to 1 when $i=j$ and equal to 0 otherwise. Then we defined
$\partial_{i} F = \partial F/\partial X^{i} = [F, P_{i}]$ and $\hat{\partial_{i}} F = \partial F/\partial P_{i} = [X^{i},F].$ In this formulation, time is not an explicit variable but the total time derivative is 
defined by another commutator with an element $H$ (the analog of the classical or quantum Hamiltonian) so that $\dot{F} = [F,H].$ Hamilton's equations are a consequence of these assumptions:
\begin{equation}
\partial H/\partial X^{i} = [H, P_{i}] = - [P_{i}, H] = - \dot{P_{i}},
\end{equation}
\begin{equation}
\partial H/ \partial P_{i} = [X^{i}, H] = \dot{X^{i}}.
\end{equation}
\noindent We then modeled dynamics by letting 
\begin{equation}
m\dot{X^i} = {\cal G}_{i} = P_{i} -  A_{i}.
\end{equation}
 where $\{ {\cal G}_{1},\cdots, {\cal G}_{d} \}$
is a collection of elements of $\cal N.$ We write ${\cal G}_{i}$
relative to the flat coordinates via ${\cal G}_{i} = P_{i} -  A_{i}.$
This is a definition of $A_{i}$ and $\partial_i F= \partial F/\partial X^{i} = [F,P_{i}].$ The formalism of gauge theory appears
naturally. In particular, defining derivations corresponding to the ${\cal G}_{i}$ by the formulas  
\begin{equation}
\nabla_{i}(F) = [F, {\cal G}_{i}],
\end{equation}
 then one has
the curvatures (as commutators of derivations) 
\begin{equation}
[\nabla_{i}, \nabla_{j}]F =[ [{\cal G}_{i}, {\cal G}_{j}],F]
\end{equation}
\begin{equation}
= [ [P_{i} -  A_{i}, P_{j} -  A_{j}],F]
\end{equation}
\begin{equation}
= [\partial_i A_j - \partial_j A_i + [A_i, A_j], F].
\end{equation}
Thus the curvature is given by the formula 
\begin{equation}
R_{ij} = \partial_{i} A_{j} - \partial_{j} A_{i} + [A_{i}, A_{j}].
\end{equation}
  We see that our curvature formula is the well-known formula for the curvature of a gauge
connection when $A_i$ is interpreted or represented as such a connection. We then saw that aspects of geometry arise in this context, including a version of the  Levi-Civita
connection.  We show how a covariant version of the Levi-Civita connection arises in this commutator calculus. This connection satisfies the formula
\begin{equation}
\Gamma_{kij} + \Gamma_{ikj} = \nabla_{j}g^{ik} = \partial_{j} g^{ik} + [g^{ik}, A_j].
\end{equation}
and so is exactly a generalization of the connection defined by Hermann Weyl in his original gauge theory \cite{Weyl}.
We compare, in Section 9, this development with the  development of gauge theory starting with Hermann Weyl. \\

A theme of this development  is the central role of the Jacobi identity 
\begin{equation}
[[X,Y],Z] + [[Y,Z],X] + [[Z,X],Y]=0
\end{equation}
in all the consequences that we draw in this noncommutative calculus. We discuss general relativity in Section 6, showing the relationship of the Bianchi identity with the Jacobi identity. In Section 7, we show how discrete calculus embeds in commutator calculus and indicate how this point of view can be used in discrete physics. \\

In all cases, studied in this paper there is the possibility for more development. By looking directing at the way calculus and physics can be done in a noncommutative world, we see that 
this sheds new light on classical mechanics, electromagnetism and gauge theory. It has been natural, since Dirac, to replace Poisson brackets by commutators and express quantum physics in 
noncommutative terms. This mode of expressing quantum mechanics can be directly accomplished using the same language as the present paper. Thus the mathematical context that we have 
expressed here is in position for interrelating classical and quantum mechanics in conceivably new ways. Some notions suggest themselves immediately such as formulating Poisson brackets directly in the noncommutative worlds. Other ideas will surely emerge as the project continues.\\

While the analogy of the Faraday tensor and the Riemann curvature tensor is clear via the commutator of the relevant derivatives, there is an asymmetry in the analogy. The coefficients of the E and M (Yang-Mills) gauge fields are effectively the difference between the generalized and geometric linear momentum $A=P-\dot{X}$ as discussed above, while the metric is the augmented background metric with momentum derivative of the gauge field A where we write (just after Lemma 3.1)
\begin{equation}
g^{ij} = [X^i, m\dot{X^j}] = [X^i, P_j] - [X^i, A_j] = \delta_{ij} - \partial A_{j}/\partial P_{i}.
\end{equation}
Such $A(P)$ dependence is to be expected \cite{G}. On other hand, the coefficients of the Levi-Civita connection, as indicated above,  are defined via the differentiation of the equation of motion. This would bring in contributions from $A(P$) via $F=dA + A\wedge A$ as part of the equation of motion. These two contributions need to be compared. Is there any room for an actual metric-based gravitational field that is not coming from $A$? See \cite{Carlip}. We are in the process of looking more closely at this context. \\
 
The present paper potentially contains  results that add to discrete approaches to the quantum gravity. Specifically by embedding discrete calculus in noncommutative calculus one makes contact with generalizations of general relativity, which includes torsion and nonmetricity in addition to the metric. 
We need to see how our work could fit with related literature such as teleparallel gravity, metric affine gravity and the geometrical formulation of quantum mechanics\cite{Hehl,Ashtekar,Rovelli,Sundance}. For example, we will see how the tetrad formalism for general relativity fits in noncommutative world context.\\

Finally, it will be important to go from the abstract algebra context in which this paper is framed to the question of physics on a noncommutative space where that space is topological.
The underlying space can be an uncountable continuum, as is traditional in classical physics, or it can be more combinatorial, but topology must be there to handle issues of connectivity and, in the quantum context, to handle issues of entanglement. There is a rich arena of questions that open from the present research.\\

We look forward to better understanding of these issues in the near future and we thank particularly David Chester and Xerxes Arsiwalla for helpful conversations.\\

\end{document}